\documentstyle[12pt]{article}

\newcommand{\const}{\mathop{\rm const}\limits}

\newcommand{\card}{\mathop{\rm card}\limits}

\newcommand{\Law}{\mathop{\rm Law}\limits}

\newcommand{\Var}{\mathop{\rm Var}\limits}

\newcommand{\supp}{\mathop{\rm supp}\limits}

\begin{document}

\begin{center}

{\bf STOCHASTIC FIELDS WITH PATHS IN ARBITRARY }

\vspace{3mm}

{\bf REARRANGEMENT INVARIANT SPACES.} \\

\vspace{3mm}

   E.Ostrovsky,  L.Sirota, {\it Israel}. \\

\vspace{2mm}

{\it Department of Mathematics and Statistics, Bar-Ilan University,
59200, Ramat Gan.}\\
e-mails: eugostrivsky@list.ru; \ sirota3@bezeqint.net \\

\vspace{3mm}
                     {\sc Abstract.}\\

\end{center}
\vspace{2mm}

 We obtain  sufficient conditions for belonging of almost all paths of  a random process   to some fixed rearrangement invariant (r.i.)
Banach functional  space, and to satisfying the Central Limit Theorem (CLT) in this  space.\par
 We describe also some possible applications.\par

\vspace{2mm}

    \par

\vspace{3mm}

\section{ Introduction. Notations. Statement of problem.}

\vspace{3mm}

 Let $  (T = \{ t \}, M, \mu) $ be measurable space with sigma - finite {\it separable} measure  $ \mu. $  Separability of the measure $ \mu $
implies separability relative a distance

$$
\rho(A_1, A_2) = \mu(A_1 \Delta A_2) = \mu(A_1 \setminus A_2) + \mu(A_2 \setminus A_1).
$$

 Let $ (L, || \cdot ||L) $ be some fixed {\it separable} Banach functional rearrangement invariant (r.i.) space over the triple
$  (T, M, \mu). $  We refer the readers about these definitions to the famous book \cite{Bennet1}, chapters 1,2.\par

 Let also $ \xi(t) = \xi(t,\omega), \ t \in T $ be {\it separable} numerical valued (real or complex) random process (r.pr.) or random field,
defined aside from the source triplet on some probability space $ (\Omega, B,{\bf P}) $  with expectation $ {\bf E } $
and variance $ \Var. $ \par

\vspace{3mm}

 {\bf We raise the question: under what the sufficient conditions almost all the paths of the r.pr. $ \xi(t)  $ belong to the space} $  L: $\par

$$
{\bf P} ( \xi(\cdot) \in L  ) = 1? \eqno(1.1)
$$

 {\bf  A second  question: under what the sufficient conditions the r.f. $  \xi(t) $ not only belongs to the space $  L $ a.e., but
 in addition satisfies the CLT in this space? } \par

\vspace{3mm}

 Recall that by definition  the r.f. $  \xi(t) $ satisfies CLT in some  Banach functional space $  B,  $ iff the suitable normed  sums

 $$
 S_n(t) := n^{-1/2} \sum_{i=1}^n \xi_i(t),
 $$
where $ \{ \xi_i(t)  \} $  are independent copies of r.f. $ \xi(t), $  converges weakly (in distribution) in the space $  B  $ as $  n \to \infty $
to a non-degenerate Gaussian random field $  S_{\infty}(t):  \ \Law(S_n) \to \Law(S_{\infty}). $  In detail:
for arbitrary continuous bounded functional $ F:  L \to R $

$$
\lim_{n \to \infty} {\bf E} F(S_n(\cdot))  = {\bf E} F(S_{\infty}). \eqno(1.2)
$$

 In particular, if $ \Law(S_n) \stackrel{def}{\to} \Law(S_{\infty}),  $ then $ \forall u = \const > 0 \ \Rightarrow $

$$
\lim_{n \to \infty} {\bf P} ( ||S_n||B > u  ) = {\bf P} ( ||S_{\infty}||B > u  ). \eqno(1.3)
$$
 Evidently, if $ \xi(t) $ satisfies the CLT in Banach functional space, then $ {\bf E} \xi(t) = 0, \ {\bf E} \xi^2(t) < \infty, \ t \in T, $
and the r.f. $ \xi(t) $ is pregaussian. This means by definition that the Gaussian centered r.f. $ S_{\infty}(t), $ which has at the same
covariation function as $ \xi(t): $

$$
R(t,s) :={\bf E} S_{\infty}(t) S_{\infty}(s)  =  {\bf E} S_{n}(t) S_{n}(s) =  {\bf E} \xi(t) \xi(s) \eqno(1.4)
$$
belongs to the space $  B $ with  probability one. \par

\vspace{3mm}

 Many sufficient conditions for the equality $ {\bf P} (\xi(\cdot) \in B) = 1 $ for different separable Banach spaces $  B  $
are obtained in \cite{Kozachenko1},  \cite{Ledoux1}, \cite{Pisier1}, \cite{Ostrovsky102}, \cite{Ostrovsky103}, \cite{Ostrovsky401},
\cite{Ostrovsky404}, \cite{Rackauskas1}, \cite{Talagrand1} - \cite{Talagrand4}.
 The case of rearrangement invariant spaces, especially ones exponential type
 Orlicz's spaces, is considered in the articles  \cite{Marcus1}, \cite{Su1}, \cite{Weber1}.\par

\vspace{3mm}

 The sufficient conditions for CLT in the Banach space of continuous functions  may be found in
\cite{Araujo1},   \cite{Dudley1}, \cite{Ledoux1}, \cite{Fortet1},
\cite{Gine1}, \cite{Gine2}, \cite{Heinkel1} etc. CLT in another separable Banach spaces is investigated, e.g. in
\cite{Billingsley2}, \cite{Garling1}, \cite{Gine2}, \cite{Ledoux1}, \cite{Pisier1}, \cite{Ostrovsky401}, \cite{Ostrovsky404},
\cite{Rackauskas1}, \cite{Zinn1}. The article \cite{Su1} is devoted to the CLT in the exponential Orlicz space, more exactly, to the
CLT in some separable subspace of the exponential type Orlicz space. \par
 The technology of application of the Banach space valued Central Limit Theorem in the parametric Monte Carlo method is described in
 \cite{Frolov1},  \cite{Grigorjeva1}, \cite{Ostrovsky1},  \cite{Ostrovsky303}. \par

\vspace{3mm}

{\sc  We need to introduce some new notations.} \par

\vspace{3mm}

{\bf A. Associate space.}  We denote by $  L' $ the {\it associate} space to the source space $ L, $ i.e. the set
of all continuous (bounded) linear functionals of the form

$$
l_g(f) = \int_T f(t) \ g(t) \ \mu(dt)
$$
with finite ordinary norm

$$
||g||L' = \sup_{f:f \in L, f \ne 0} \left[ \frac{l_g(f)}{||f||} \right] =
 \sup_{f:f \in L, ||f|| = 1} l_g(f). \eqno(1.5)
$$
 It is known that $ (L')' = L,  $ see \cite{Bennet1}, chapter 2. \par

 Denote by $ S = S_e $  the set of all {\it extremal} points of the unit surface of associate space $  L', $ so that $ S = S_e $
is symmetric: $ - S = S $ and  $ \forall g \in S \ \Rightarrow ||g||L' = 1; $  then

$$
\forall f \in L \ \Rightarrow ||f||L = \sup_{g \in S_e} l_g(f). \eqno(1.6)
$$

 \vspace{3mm}

{\bf B. Fundamental function.}  Recall that the fundamental function $ \phi(L, \delta), \ \delta = \const > 0 $ for r.i. Banach
functional  space $ (L, ||\cdot||L) $ is defined as follows. Let $  A $ be some measurable set with measure
$ \delta: \ \mu(A) = \delta > 0.  $ Then

$$
\phi(L, \delta) \stackrel{def}{=} ||I(A)||L; \eqno(1.7)
$$
in the sequel $ I(A) = I_A(t) $ denotes an indicator function of the set $  A. $ This definition does not dependent on the concrete
representation of the set $  A $ and play a very important role in the theory of operators, theory of Fourier series etc.see \cite{Bennet1},
chapters 4,5. \par

 Note that $ \phi(L, \delta) \cdot \phi(L', \delta) = \delta.  $\par

\vspace{3mm}

{\bf C. Metric entropy.}

 Further, let  $ (X, \rho) $ be compact metric space relative the distance (or semi - distance) function  $ r = r(x_1, x_2). $
 Denote as usually by $ N(X,r,\epsilon), \ \epsilon = \const > 0 $ the minimal number of closed $  r \ -  $ balls of a radii
 $ \epsilon $ which cover all the set $  X; $ obviously, $ \forall \epsilon > 0 \ \Rightarrow N(X,r,\epsilon) < \infty. $ \par
  The quality

 $$
 H(X,r,\epsilon) = \ln N(X,r,\epsilon) \eqno(1.8)
 $$
is named {\it entropy  } of the set $  X  $ relative the distance $  r  $ at the point $ \epsilon. $ \par
 This notion is in detail investigated, e.g., in \cite{Ostrovsky1},  chapter 3, section 3.2. \par

\vspace{3mm}

{\bf D. Grand Lebesgue spaces.  }

  Recently, see \cite{Fiorenza1}, \cite{Fiorenza2}, \cite{Fiorenza3}, \cite{Iwaniec1}, \cite{Iwaniec2}, \cite{Kozachenko1},
 \cite{Ostrovsky1}  etc.  appears the so-called {\it Grand Lebesgue Spaces} $ GLS = G(\psi) =  G(\psi; a,b), \ a,b = \const,
 a \ge 1, a < b \le \infty, $ spaces consisting
 on all the measurable functions $ f: T \to R, $ where (recall) $ (T =  \{t \}, M, \mu) $ is measurable space with non - trivial
 sigma - finite measure $  \mu, $ having a finite norms

$$
  ||f||G(\psi) = ||f||G(\psi; \mu)  \stackrel{def}{=} \sup_{p \in (a,b)} \left[ |f|_{p,\mu} /\psi(p) \right], \eqno(1.9)
$$
where we define as usually

$$
|f|_{p,\mu}:= \left[ \int_T |f(t)|^p \ \mu(dt)  \right]^{1/p}, \ 1 \le p < \infty.
$$
and we define correspondingly for a random variable $ \eta $

$$
|\eta|_p := \left[ {\bf E} |\eta|^p  \right]^{1/p}.
$$

 Here $ \psi(\cdot) $ is some continuous positive on the {\it open} interval  $ (a,b) $ function such that

$$
     \inf_{p \in (a,b)} \psi(p) > 0.
$$

We can suppose without loss of generality

$$
\inf_{p \in (a,b)} \psi(p) = 1.
$$

 Notation: $  (a,b) = \supp \psi. $\par

 As the capacity of the measure $ \mu $ may be picked the probability measure $ {\bf P} $ defined on all  the
measurable sets $  B; \ B \subset \Omega. $ \par

  This spaces are rearrangement invariant, see \cite{Bennet1}, and
are used, for example, in the theory of probability \cite{Talagrand1}, \cite{Kozachenko1}, \cite{Ostrovsky1};
theory of Partial Differential Equations \cite{Fiorenza2}, \cite{Iwaniec2}; functional analysis; theory of Fourier
series; theory of martingales  etc.\par

 Let $ \delta = \const > 0; $ the fundamental function  $ \phi(\delta) = \phi_{ G(\psi)}(\delta)  $  of the space $  G(\psi) $
may be calculated as follows:

$$
\phi_{ G(\psi)}(\delta) = \sup_{p \in (a,b)} \left[ \frac{\delta^{1/p}}{\psi(p)} \right].
$$

 The fundamental function  for GLS is in detail investigated in \cite{Holon1}; in particular, it was therein  calculated many examples.
Roughly speaking, there is an essential difference for behavior of this function between the cases $  b < \infty $  and $ b = \infty. $\par

 \vspace{3mm}

 {\bf D. \ Definition 1.1.}  Suppose the set $ T $ is equipped in addition with some distance (semi-distance) function $ d = d(t,s) $ such that
 sigma field $  M $ is Borelian sigma algebra (metric measurable space) and such that the metric space $  (T,d) $ is compact space. \par

 The {\it  ball function  }  $  r(\delta) = r(T, \delta), \ \delta = \const > 0 $ for the measurable space $ (T,M,\mu)  $
equipped  with a distance  $ d = d(t,s) $ is by definition the function of a form

$$
r(T, d, \delta)= r(T, \delta) = r(\delta) \stackrel{def}{=} \sup_{t \in T} \mu  [B(t, \delta)], \eqno(1.10)
$$
where as ordinary

$$
B(t,\delta) = B(t,d, \delta) = \{ z, z \in T, \ d(t,z) \le \delta \}
$$
is closed $ \delta \ -  $ ball in the metric space $  (T,d)  $ with the center at the point $  t, \ t \in T. $ \par

 It is clear that $ \lim_{\delta \to 0+} r(T,\delta) = 0.  $

\vspace{3mm}

 Suppose now the measure $ \mu $ in  is probabilistic: $ \mu(T) = 1, $
and let the function $ \psi = \psi(p) $  be such that $ \supp \psi = (1,\infty). $ Then the  Grand Lebesgue Space $ G(\psi) $
coincides up to norm equivalence to the subspace of all mean zero:  $ \int_T f(t) \ \mu(dt) = 0  $
 measurable function (random variables) of the so-called {\it exponential Orlicz space  } $ L(N) = L(N; T,\nu) $
 with {\it  exponential } Orlicz-Young function $ N = N(u),$ and conversely proposition is also true: arbitrary exponential
Orlicz space $ L(N) $  coincides with some Grand Lebesgue Space, see \cite{Kozachenko1}.\par

 In detail, introduce the function $ \phi(\cdot) $ as follows:

$$
\chi(p) := p \ \psi(p), \ p \ge 2, \ \chi(p) = C \cdot p^2, \ 0 \le p \le 2, \ C: 2 \chi(2) = 4 C,
$$

$$
\phi(y) :=  (\chi(y))^{-1}, \ y \ge 0; \ \phi(y) := \phi(|y|), \ y < 0;
$$
then

$$
N(u) = \exp( \phi^*(u) ) - 1,
$$
where $ \phi^*(\cdot) $ is Young - Fenchel, or Legendre  transform for the function $ \phi: $

$$
\phi^*(u)  = \sup_{y \in R} (uy - \phi(y)),
$$
see \cite{Ostrovsky1}, chapter 1, theorem 1.5.1.\par

 \vspace{3mm}

 The finiteness of the $ G\psi $ norm for the r.v. $  \xi $  allows to obtain the quasy - exponential bounds for its tail of distribution.
Indeed, if we denote

$$
T_{\xi}(x) = \max \left[ {\bf P}(\xi > x), \    {\bf P}(\xi < x)  \right], \ x \ge 2,
$$
and if $  0 < ||\xi||G\psi < \infty, $ then

$$
T_{\xi}(x) \le \exp \left\{ -  [ p \ln \psi(p) ]^*(\ln x/||\xi||G\psi) \right\},  \ x \ge 2 ||\xi||G\psi,
$$
and the conversely proposition is true. Namely, if for some r.v. $  \xi  $

$$
T_{\xi}(x) \le \exp \{ - h(\ln x)    \}, \ x \ge 2,
$$
where $ h = h(y) $ is positive continuous convex strictly monotonically  increasing function such that
$ \lim_{y \to \infty} h(y) = \infty, $ then

$$
||\xi||G\psi \le C(h) < \infty, \ {\bf where} \  \psi(p) := \exp \left( \frac{h^*(p)}{p}  \right).
$$

 Let for instance $ m = \const > 0 $ and define  the following $  \psi $ function

$$
\psi_m(p) = p^{1/m}, \ 1 \le p < \infty.
$$
 The r.v. $ \psi $ belongs to the space $  G\psi_m $ iff for some positive constant $  C(m) $

$$
T_{\xi}(x) \le \exp  \left( - C(m) \ x^m   \right), \ x \ge 2.
$$

 The last proposition is well known, see e.g. \cite{Su1}. \par

The case $  m = 2 $  correspondent to the case of the so-called subgaussian random variables, centered or not. Here
$ \psi(p) = \sqrt{p}. $ \par

\vspace{3mm}

 Let now $ F = \{  f_{\alpha}(t) \}, \ \alpha \in A $ be a {\it family} of measurable functions such that

 $$
 \exists (a,b), \ 1 \le a < b \le \infty, \ \ \forall p \in (a,b) \ \Rightarrow \sup_{\alpha \in (a,b)} ||f_{\alpha}||_{p, \mu} < \infty. \eqno(1.11)
 $$
 The function
$$
\psi_F(p) := \sup_{\alpha \in (a,b)} || \ f_{\alpha}(\cdot) \ ||_{p, \mu} \eqno(1.12)
$$
 is said to be {\it natural function  } for the family $  F. $ This function is obviously minimal up to equivalence function $ \psi $ for which

$$
\sup_{\alpha \in A}  ||f_{\alpha}||G \psi = 1.
$$

\vspace{3mm}

\section{ First condition.}

\vspace{3mm}

 {\it We do not assume in this and in the next sections that the r.i. space $  L  $ is separable. } \par

 \vspace{3mm}

  Suppose here that the r.f. $ \xi(t) $ belongs uniformly in $ t, \ t \in T $ to some non-trivial  $  G\psi_0  $ space:

$$
 \exists a,b: 1 = a  < b \le \infty, \ \Rightarrow \forall p \in (a,b) \  \psi_0(p) := \sup_{t \in T} ||\xi(t)||_{L(p), \Omega} < \infty. \eqno(2.1)
$$

 In what follow we can use instead the natural function $ \psi_0 $ in (2.1) arbitrary its majorant $ \psi = \psi(p)  $ from
the set $  G \Psi $  with at the same support $ (a,b). $ \par

 Let us introduce a so-called {\it natural, } i.e. generated  only by means of the values of the r.f. $ \{\xi(t) \}, \ t \in T, $
  on the set $  T  $ {\it bounded } semi-distance  $  d_{\psi} = d_{\psi}(t,s) $ as follows

$$
d_{\psi}(t,s) = d_{\psi} \stackrel{def}{=} ||\xi(t) - \xi(s)|| G\psi. \eqno(2.2)
$$
 Note that for natural distance $  d_{\psi_0}  \ \hspace{6mm} d_{\psi_0}(t,s) \le 2 $ and that the r.f. $ \xi(t) $ is stochastic continuous
relative this distance.\par

\vspace{3mm}

{\bf Theorem 2.1.}  {\it Suppose that for some $ q = \const \in (0,1) $ the following entropy series converge: }

$$
\sigma = \sigma(q) \stackrel{def}{=}   \sum_{n=0}^{\infty} q^{n} \ N \left(T, d_{\psi}, q^{n+1} \right) \  r \left(T, q^{n} \right) < \infty. \eqno(2.3)
$$

 {\it  Then  }

 $$
 {\bf P} (\xi(\cdot) \in L) = 1 \eqno(2.4)
 $$

 {\it and moreover if in addition $ \sup_{ t \in T } ||\xi(t)||G\psi = 1,  $ then   }

$$
|| \ ||\xi(\cdot)||L \ ||G\psi \le  \underline{\sigma} := \inf_{q \in (0,1)}
\left\{ \sum_{n=0}^{\infty} q^{n} \ N \left(T, d_{\psi}, q^{n+1} \right) \ r \left(T, q^{n} \right) \right\}. \eqno(2.5)
$$

\vspace{3mm}

{\bf Proof.}  We can and will assume without loss of generality $ \sigma = 1. $ Further, denote by $  T_n = T_n \left(q^{n} \right)  $
the minimal $  q^{n} $ of the set $ T $ relative the distance $  d_{\psi}. $ \par
 This net in not necessary to be unique; we pick arbitrary fixed but non-random one. \par
We have  on the basis of entropy definition

$$
\card(T_n) = N(q^{n+1} ) := N(T, d_{\psi}, q^{n+1}). \eqno(2.6)
$$
We define for arbitrary element $ t \in T  $ and any value $  n = 0,1,2, \ldots  $ the following "projection" $ \theta_n(t): $

$$
d_{\psi}(t, \theta_n(t)) \le q^{n}, \ \theta_n(t) \in T_n. \eqno(2.7)
$$
 This point may be also not unique, but we choose it non - random.  By definition, $ \theta_0(t) := t_0 \in T_0 $
 be some fixed point inside the set $  T.  $ \par
  We have

 $$
 \xi(t) = \sum_{n=0}^{\infty} \left( \xi(\theta_{n+1}(t) - \xi(\theta_n(t)  \right), \eqno(2.8)
 $$
therefore

$$
||\xi(\cdot)||L \le \sum_{n=0}^{\infty} \eta_n, \ \eta_n := ||  \xi(\theta_{n+1}(\cdot)) - \xi(\theta_n(\cdot)) ||L. \eqno(2.9)
$$

 The function $ t \to \xi(\theta_{n+1}(\cdot)) - \xi(\theta_n(\cdot))  $ is simple (stepwise),  therefore it belongs to the  space $  L. $
 The amount of the $  d_{\psi} $
balls of a radii $ q^{n+1} $ is less or equal than $ N(q^{n+1}). $ The value $ |\xi(\theta_{n+1}(\cdot)) - \xi(\theta_n(\cdot))| $
does not exceed   the value $ q^{n}. $\par
 Since the $ L \ - $ space is rearrangement invariant,

$$
||\eta_n|| G\psi \le q^n \ N(q^{n+1}) \ r(T, q^n). \eqno(2.10)
$$
 It remains to use the triangle inequality and completeness of the r.i. Banach functional space $ L:  $\par

$$
|| \ || \xi(\cdot)||L \ ||G\psi \le \sum_{n=0}^{\infty}  q^n \ N(q^{n+1}) \ r(T, q^n) = \sigma(p).
$$
 Since the value $ q $ is arbitrary inside the interval $ (0,1), $

 $$
{\bf P} ( \xi(\cdot) \in L) = 1, \ \hspace{6mm} || \ || \xi(\cdot)||L \ ||G\psi \le \inf_{q \in (0,1)} \sigma(q) = \underline{\sigma},
 $$
Q.E.D. \par

 As a consequence: under formulated above conditions $ {\bf P} ( \xi(\cdot) \in L) = 1 $ and moreover

$$
{\bf P} ( || \xi(\cdot)||L > x )
 \le \exp \left\{ -  [ p \ln \psi(p) ]^*(\ln x/\underline{\sigma}) \right\},  \ x \ge 2 \underline{\sigma}. \eqno(2.11)
$$

\vspace{3mm}

{\bf  Remark 2.1.} The expression  $ || \ || \xi(\cdot)||L \ ||G\psi  $ is called mixed, or on the other words Bochner's
norm for the function of two variables (random process) $ \xi = \xi(t, \omega). $ \par

\vspace{3mm}

\begin{center}

{\bf  Examples.}

\end{center}

\vspace{3mm}

 {\bf First example.} \par

\vspace{3mm}

   Conditions.  Let $  \psi = \psi(p)  $ be some non-trivial: $ b =\sup \{p; \ p \in \supp \psi\} > 1  $ natural
 function for the r.f. $ \xi(t): \ \sup_t ||\xi(t)||G\psi = 1. $ Then evidently $  d(t,s) = ||\xi(t) - \xi(s)||. $  \par
  Suppose first of all that $  r(T,\delta) \le \delta^s, \delta \in (0,1)  $ for some positive value $ s = \const  > 0. $ \par
If  for instance $ T $  is closure of an open set in the space $  R^d $ and $  r(t,s) \asymp |t-s|^{\alpha}, $ where  $|t|$ is usually
Euclidean norm  and $ \alpha = \const \in (0,1], $ then $  s = d/\alpha. $\par
 Assume further

$$
N(T,d_{\psi}, \epsilon) \le \epsilon^{-\kappa}, \ \epsilon \in (0,1), \ \kappa = \const \in (0, 1 + s). \eqno(2.12)
$$
  The last equality (2.12) implies that the entropy dimension of the set $  d  $ relative the distance $ d_{\psi} $ is restricted

$$
\dim_{d_{\psi} }(T)   = \kappa < 1 +s.
$$
 In the sequel we assume the values $ \kappa,s $ to be fixed.  Another notations: $ \Delta = 1 + s - \kappa > 0, \ \lambda = \kappa /\Delta. $ \par

 We deduce that all the conditions of theorem 2.1 are satisfied and we obtain
 after some calculations that the optimal value of the parameter $ p $ is following:

$$
p_0 = \left(  \frac{\lambda}{1 + \lambda}  \right)^{1/\Delta} = \left[ \frac{\kappa}{\kappa + \Delta} \right]^{1/\Delta}
$$
and correspondingly

$$
\underline{\sigma} = \left[\frac{\kappa}{\Delta}\right]^{- \kappa/\Delta  } \cdot
\left[ 1 + \frac{\kappa}{\Delta} \right]^{ - 1 - \kappa/\Delta }. \eqno(2.13)
$$

\vspace{3mm}

 {\bf Second example.} \par

\vspace{3mm}

 Al the parameters are as before in the first example aside from the entropy condition:

$$
N(\epsilon) \le \epsilon^{-(1+s)} \ |\ln \epsilon|^{-\beta}, \ \epsilon \in (0, \ 1/e), \ \beta = \const > 1, \ \kappa = 1 + s. \eqno(2.14)
$$

 We deduce again that all the conditions of theorem 2.1 are satisfied and we obtain
 after some calculations that the optimal value of the parameter $ p $ is following:

$$
p_0 = \exp (-\beta/\kappa)
$$
and correspondingly

$$
\underline{\sigma}  \le e^{-\beta} \ \beta^{\beta} \ \kappa^{-\beta} \ (\zeta_R(\beta) -1), \eqno(2.15)
$$
where $ \zeta_R(\cdot) $ denotes ordinary Rieman's zeta function.\par

\vspace{3mm}

\section{ Second condition.  }

\vspace{3mm}

   Since the r.f. $  \xi(t) $ is stochastic continuous, it is $ (T,M) $ measurable. \par

\vspace{3mm}

  {\it We suppose in addition that for arbitrary (non-random!)
function $ g = g(t) $  from the space $ L' $ there exists (with probability one) the following linear functional (integral):}

$$
\forall g \in L' \hspace{4mm}  \exists l_{\xi}(g):= \int_T \xi(t) \ g(t) \ \mu(dt).\eqno(3.1)
$$
  A simple  sufficient condition for (3.1)  is following: the function $ t \to {\bf E} |\xi(t)| $ there exists and belongs to
the space $  L.$ \par
 It remains to establish the finiteness with probability one the value

$$
\nu := \sup_{g: ||g||L' = 1} l_{\xi}(g).
$$
 As long as the set $  S = S_e $ is the set of all extremal points in the centered unit ball of the space $  L', $

$$
\nu = \sup_{g \in S} l_{\xi}(g). \eqno(3.2)
$$
{\it   In what follow we consider in this section only the case  when} $ g(\cdot) \in S. $\par

 Suppose that as in the last section that the family of random variables  $ \{  \xi(t) \}, \ t \in T $ obeys some non - trivial
natural $  \psi $ function:

$$
\sup_{t \in T} || \xi(t)||G\psi = 1, \ b := \sup \supp \psi > 1. \eqno(3.3)
$$

 We estimate using triangle (Marcinkiewicz) inequality

$$
|| l_{\xi}(g)||G\psi  = || \int_T \xi(t) \ g(t) \ \mu(dt)  ||G\psi \le
 \int_T || \xi(t)|| G\psi \ |g(t)| \ \mu(dt)   \le
$$

$$
\int_T |g(t)| \mu(dt) = ||g(\cdot)||L_1(T,\mu) \le C_1 = \const < \infty.\eqno(3.4)
$$

We introduce now the (semi-) distance $ \rho(g_1,g_2) $ on the set $  S_e $ as follows:

$$
\rho(g_1,g_2) := || g_1(\cdot) - g_2(\cdot)||L_1(T, \mu).
$$
 We have analogously

$$
|| l_{\xi}(g_1)  -  l_{\xi}(g_2) ||G\psi  = || \int_T \xi(t) \ (g_1(t) - g_2(t)) \ \mu(dt)  ||G\psi \le
$$

$$
 \int_T || \xi(t)|| G\psi \ |g_1(t) - g_2(t)| \ \mu(dt)   = ||g_1 - g_2||L_1(T, \mu) = \rho(g_1, g_2). \eqno(3.5)
$$

\vspace{3mm}

 Define also for arbitrary function $ f: R_+ \to R $ the Young-Fenchel co-transform $ f_* $ by an equality

$$
f_*(x) \stackrel{def}{=} \inf_{y \ge 0} (xy + f(y)), \eqno(3.6)
$$
and introduce the diameter of the set $  S  $ relative the semi-distance $ \rho(\cdot, \cdot) $
$ D := \sup_{g_1,g_2 \in S} \rho(g_1,g_2) < \infty $ and a function  $ v(y) = \ln \psi(1/y).  $\par

\vspace{3mm}

{\bf Theorem 3.1.}  {\it  Suppose in addition to the formulate above conditions}

$$
I := \int_0^D \exp( v_*(2 + \ln N(T, d_{\psi}, \epsilon ))) \ d \epsilon < \infty. \eqno(3.7)
$$
{\it Then }  $ {\bf P}(\xi(\cdot) \in L) = 1.  $ \par

\vspace{3mm}

{\bf Proof \ } follows immediately from the theorem 3.17.1 of  the monograph \cite{Ostrovsky1}, chapter 3, section 3.17, where is proved
in particular that

$$
(|| \ ||\xi||L \ ||G\psi = ) \hspace{4mm} ||\nu||G\psi = || \sup_{g \in S} l_{\xi}(g) ||G\psi \le 9 I < \infty. \eqno(3.8)
$$
 This completes the proof of theorem 3.1.\par

\vspace{3mm}

\section{ Conditions for the Central Limit Theorem.}

\vspace{3mm}

 {\it We suppose in this section  that the r.i. space $  L  $ is separable. } \par

\vspace{3mm}

{\it Suppose in addition that the random field $  \xi(t) $ is mean zero,  has uniform bounded  second moment and is pregaussian.} \par

\vspace{3mm}

{\bf First version.}\par

\vspace{3mm}

  Suppose here in this subsection as before that the r.f. $ \xi(t) $ belongs uniformly in $ t, \ t \in T $ to some non - trivial  $  G\psi_0  $ space:

$$
 \exists a,b: 2 = a  < b \le \infty, \ \Rightarrow \forall p \in (a,b) \  \psi_0(p) := \sup_{t \in T} ||\xi(t)||_{L(p), \Omega} < \infty. \eqno(4.1)
$$

 In what follow we can use instead the natural function $ \psi_0 $ in (4.1) arbitrary its majorant $ \psi = \psi(p)  $ from
the set $  G \Psi $  with at the same support $ (a,b). $ \par

  We define for arbitrary such a function $ \psi(\cdot) $ its {\it Rosenthal's}  transform $  \psi_R(\cdot): $

$$
\psi_R(p) \stackrel{def}{=} \frac{C_R \ p}{ \ln p} \cdot \psi(p), \ p \in (2,b), \ C_R := 1.77638.\eqno(4.2)
$$

 It is clear that if $  b <  \infty, $ then $  \psi_R(\cdot) \asymp \psi(\cdot), \ 1 \le p < b.   $ Therefore, we will assume in this approach
 $  b = \infty. $ \par

 The classical Rosenthal's inequality \cite{Rosenthal1} asserts in particular that if $  \{ \zeta_i \}, \ i = 1,2,\ldots $ are
the sequence of i., i.d. {\it centered} r.v. with finite $  p^{th} $ moment, then

$$
\sup_n \left| n^{-1/2} \sum_{i=1}^n \zeta_i  \right|_p  \le \frac{C_R \ p}{ \ln p} \ |\zeta_1|_p, \ p \ge 2. \eqno(4.3)
$$
 About the exact value of the constant $  C_R $ see the article \cite{Ostrovsky601}.  Note that for symmetrical distributed r.v.
 $  C_R \le 1.53573. $\par

 \vspace{3mm}

 Let us consider the normed sums

$$
 S_n(t) := n^{-1/2} \sum_{i=1}^n \xi_i(t), \ n = 1,2,\ldots. \eqno(4.4)
$$

 It follows from Rosenthal's inequality

$$
\sup_n \sup_{t \in T} || S_n(t) ||G\psi_R \le 1,  \eqno(4.5a)
$$

and we define

$$
\rho_{\psi}(t,s) \stackrel{def}{=}  \sup_n || S_n(t) - S_n(s) ||G\psi_R.  \eqno(4.5b)
$$

\vspace{3mm}

{\bf Theorem 4.1.}  {\it Suppose that for some $ q = \const \in (0,1) $ the following entropy series converge: }

$$
\gamma = \gamma(q) \stackrel{def}{=}   \sum_{n=0}^{\infty} q^{n} \ N \left(T, \rho_{\psi}, q^{n+1} \right) \  r \left(T, \rho_{\psi}, q^{n} \right) < \infty. \eqno(4.6)
$$

 {\it  Then  the sequence of r.f $ \xi_i(t), \ t \in T  $ satisfies the CLT in the r.i. $  L  $ space. }\par

 \vspace{3mm}

 {\bf Proof.} The convergence of finite dimensional (cylindrical)  distributions of r.f. $  S_n(t) $ to the finite dimensional of the Gaussian r.f.
$ S_{\infty}(t), $  which sample path in turn belongs to the space $  L  $ is evident. It remains to establish the weak compactness of
measures in the space $  L  $ generated by r.f. $  S_n(\cdot). $\par

 We apply theorem 2.1  to the random field $  S_n(\cdot) $ for arbitrary fixed  value $  q  $ for which $ \gamma(q) < \infty: $

$$
\sup_n || \ ||S_n(\cdot)||L \ ||G\psi_R \le  \gamma(q) < \infty. \eqno(4.7)
$$

  Further, as long as the space $  L  $ is presumed to be separable, there exists a {\it compact} linear operator $ U: L \to L $ such that
$  U^{-1}\xi \in L, \ \Rightarrow U^{-1} S_n \in L   $ and moreover

$$
\sup_n || \ ||U^{-1} S_n(\cdot)||L \ ||G\psi_R \le  1,  \eqno(4.8)
$$
see \cite{Ostrovsky304}, \cite{Ostrovsky602}. Therefore,

$$
\lim_{Z \to \infty}  \sup_n {\bf P} \left( || U^{-1} S_n ||L > Z  \right) = 0. \eqno(4.9)
$$
 Since  the set

$$
 W(Z) = \{ f, \ f \in L, \ \ ||U^{-1}f||L \le Z  \}
$$
is compact subset  of the space $  L, $ the equality (4.9) proves the proposition of theorem 4.1.\par
See also the criterion for the functional CLT in the famous book  \cite{Ledoux1}, chapter 6. \par

\vspace{3mm}

{\bf Second version.}\par

\vspace{3mm}

 We intent to mention and to generalise the third our section.

\vspace{3mm}

   {\it We suppose as before  in addition that for arbitrary (non - random!)
function $ g = g(t) $  from the space $ L' $ there exists (with probability one) the following linear functional (integral):}

$$
\forall g \in L' \hspace{4mm}  \exists l_{\xi}(g):= \int_T \xi(t) \ g(t) \ \mu(dt).
$$
  Then automatically

$$
\forall g \in L' \hspace{4mm}  \exists l_{S_n}(g):= \int_T S_n(t) \ g(t) \ \mu(dt).
$$

  It remains only to establish the finiteness with probability one the value

$$
\lambda := \sup_n \sup_{g: ||g||L' = 1} l_{S_n}(g).
$$
or equally

$$
\lambda = \sup_{g \in S} l_{S_n}(g).
$$

{\it   In what follow we consider in this section only the case  when} $ g(\cdot) \in S. $\par

 Suppose that as in the third section that the family of random variables  $ \{  \xi(t) \}, \ t \in T $ obeys some non-trivial
natural $  \psi $ function:

$$
\sup_{t \in T} || \xi(t)||G\psi = 1, \ b := \sup \supp \psi > 1.
$$
 Then

$$
\sup_{t \in T} || S_n(t)||G\psi_R = 1, \ b := \sup \supp \psi  = \infty,
$$
the case $  b < \infty $ is trivial. \par

 We estimate using triangle (Marcinkiewicz) inequality

$$
|| l_{S_n}(g)||G\psi_R  = || \int_T S_n(t) \ g(t) \ \mu(dt)  ||G\psi_R \le
 \int_T || S_n(t)|| G\psi_R \ |g(t)| \ \mu(dt)   \le
$$

$$
\int_T |g(t)| \mu(dt) = ||g(\cdot)||L_1(T,\mu) \le C_1 = \const < \infty.\eqno(4.10)
$$

 Recall that we introduced  the (semi-) distance $ \rho(g_1,g_2) $ on the set $ S = S_e $ as follows:

$$
\rho(g_1,g_2) := || g_1(\cdot) - g_2(\cdot)||L_1(T, \mu).
$$
 We have analogously

$$
|| l_{S_n}(g_1)  -  l_{S_n}(g_2) ||G\psi_R  = || \int_T S_n(t) \ (g_1(t) - g_2(t)) \ \mu(dt)  ||G\psi_R \le
$$

$$
 \int_T || S_n(t)|| G\psi_R \ |g_1(t) - g_2(t)| \ \mu(dt)   = ||g_1 - g_2||L_1(T, \mu) = \rho(g_1, g_2). \eqno(4.11)
$$
 and introduce the diameter of the set $  S  $ relative the semi-distance $ \rho(\cdot, \cdot) $
$ D := \sup_{g_1,g_2 \in S} \rho(g_1,g_2) < \infty $ and a function  $ v_R(y) = \ln \psi_R(1/y).  $\par

\vspace{3mm}

{\bf Theorem 4.2.}  {\it  Suppose in addition to the formulate above conditions}

$$
J := \int_0^D \exp( v_{R,*}(2 + \ln N(T, d_{\psi_R}, \epsilon ))) \ d \epsilon < \infty. \eqno(4.12)
$$
{\it Then the sequence of the r.f. $ S_n(\cdot) $ satisfies the CLT in the space } $  L. $   \par

\vspace{3mm}

{\bf Proof \ } follows immediately from the theorem 3.17.1 of  the monograph \cite{Ostrovsky1}, chapter 3, section 3.17, where it is proved
in particular that uniformly in $  n  $

$$
( \sup_n || \ ||S_n||L \ ||G\psi_R = ) \hspace{4mm} || \lambda||G\psi_R = \sup_n || \sup_{g \in S} l_{S_n}(g) ||G\psi_R \le 9 J < \infty.
$$
 This completes the proof of theorem 4.2.\par

\vspace{3mm}

\section{ Concluding remarks. Applications.}

\vspace{3mm}

{\bf A. Applications in the Monte-Carlo method. } \par

\vspace{3mm}

 Let us consider here the problem of Monte-Carlo approximation and construction
of a confidence region in the  $ L \ - $ space norm for the parametric integral of a view

$$
I(t) = \int_X g(t,x) \ \nu(dx). \eqno(5.1)
$$

 Here $ (X, F, \nu) $ is also a probabilistic space with normed: $  \nu(X) = 1 $ non-trivial
measure $ \nu. $ \par

 A so-called ”Depending Trial Method” estimation for the integral (5.1) was introduced by Frolov A.S.and Tchentsov N.N., see
 \cite{Frolov1}:

$$
I_n(t) = n^{-1} \sum_{i=1}^{\infty} g(t, \eta_i), \eqno(5.2)
$$
 where $  \{ \eta_i \} $ is the sequence of $  \nu $ distributed: $ {\bf P} (\eta_i \in A) = \nu(A) $
independent random variables.

 Suppose that the sequence of r.f. $ g(t, \eta_i) - I(t)  $ satisfies the CLT in some Banach r.i. space $  L; $ then

$$
\lim_{n \to \infty} {\bf P} \left( \sqrt{n} || I_n(\cdot) - I(t)||L > u \right) = {\bf P} ( ||\zeta(\cdot)||L > u), \ u > 0; \eqno(5.3)
$$
therefore

$$
 {\bf P} \left( \sqrt{n} || I_n(\cdot) - I(t)||L > u \right) \approx {\bf P} ( ||\zeta(\cdot)||L > u), \ u > 0; \eqno(5.4)
$$
 The last equality may be used by the construction  of a confidence region  (c.r.) in the $  L  $ norm for the integral $  I(t).$ Namely,
equating the right - hand side  of (5.3a) to some "small" number $ \delta, $  for instance $ \delta = 0.05 $  or $ \delta = 0.01 $ etc., where
the value $  1 - \delta  $ is reliability of the c.r.:

$$
{\bf P} ( ||\zeta(\cdot)||L > u_0) = \delta, \eqno(5.5)
$$
we obtain an asymptotical c. r. of a form: with probability  $  \approx 1 - \delta  $

$$
|| I_n(\cdot) - I(t)||L \le  \frac{u_0}{\sqrt{n}}. \eqno(5.6)
$$
 See for detail description the articles  \cite{Grigorjeva1}, \cite{Ostrovsky302}, \cite{Ostrovsky104}, \cite{Ostrovsky105}.\par

\vspace{3mm}

{\bf B.  A case of H\"older - Lipshitz space.} \par

\vspace{3mm}

 The CLT in the so-called H\"older (Lipshitz) space  $ H^o(\omega),  $ (but which is not rearrangement invariant),
 is investigated, e.g. in \cite{Ostrovsky1}, chapter 4, section 4.13. \par
Recall that the  H\"older {Lipshitz} space  $ H^o(\omega)  $  consists on all the numerical continuous relative some distance $ d = d(t,s) $
functions $ f: T \to R $  satisfying the condition

$$
\lim_{\delta \to 0+} \frac{\omega(f, \delta)}{\omega(\delta)} = 0.  \eqno(5.7)
$$
 Here $ \omega(f, \delta)  $ is uniform module of continuity of the function $  f:$

$$
\omega(f, \delta) = \sup_{t,s: d(t,s) \le \delta} |f(t) - f(s)|,
$$
$ \omega(\delta)  $ is {\it  some} non - trivial (continuous) module of continuity. For example, $ \omega(0+) = \omega(0) = 0,
\ \delta > 0 \ \Rightarrow \omega(\delta) > 0 $ etc.\par
 The metric space $ (T,d) $ is presumed to be compact. \par

The norm of the space $ H^o(\omega) $  is defined as follows:

$$
||f||H^o(\omega) = \sup_{t \in T} |f(t)| + \sup_{\delta \in (0,1)} \omega(f, \delta). \eqno(5.8)
$$

 This modification of the classical Lipshitz space is in general case separable. \par
The recent version for CLT in H\"older spaces, for example for the Banach space valued random processes, see in
\cite{Rackauskas1}. \par

 In the article of  B.Heinkel \cite{Heinkel1}  is obtained sufficient condition for  CLT in the space of continuous functions
$  C(T,d) $  in the more modern terms of "majorizing measures";  see \cite{Fernique1}, \cite{Talagrand1} - \cite{Talagrand4}. \par

 It is interest by our opinion to obtain the conditions for CLT in these terms for the H\"older-Lipshitz spaces, as well as for the
separable functional rearrangement invariant spaces. \par

\vspace{3mm}

{\bf C. Counterexample.}

\vspace{3mm}

  Let $  T = [0, \ 2 \pi].  $ There exists an example of mean zero {\it continuous  } periodical r.pr. $  \xi(t) $  constructed by
means of lacunar trigonometrical series   which does not satisfy the CLT in the space  $ C(T),  $ see \cite{Kozachenko1}. The analog
of the conditions of theorem 3.1 for this space is satisfied but the conditions of theorem 4.2 are not.\par
  This process can serve as an example (counterexample) to our situation. More detail, let us consider the Orlicz space $ \hat{L} $
over the set $  T  $ with the Young-Orlicz function $  N(u) = \exp(u^4) - 1. $ As we know, the norm in this space may be defined
up to equivalence as follows:

$$
||f|| \hat{L} := \sup_{p \ge 1} \left[ \frac{|f|_p}{p^{ 1/4 }} \right].
$$
 But this space is not separable. In order to obtain the separable space, we introduce as a capacity of the space $  L  $ the subspace
of $ \hat{L}  $  consisting on all the function $ f \in \hat{L} $ for which

$$
\lim_{p \to \infty} \left[ \frac{|f|_p}{p^{ 1/4 }}  \right] = 0.
$$
 As  long as the limiting Gaussian process $  S_{\infty}(t), \ t \in T $ described in third section does not belongs to the space $  L, $ the
 continuous a.e. r.pr. $ \xi(t) $  does not satisfy the CLT also in the space  $ L.  $ \par

\end{document}